A Mathematical Problem-Solving Pipeline (MPSP) to strengthen scaffolding in higher education STEM courses

Rochelle E. Tractenberg, ORCID: 0000-0002-1121-2119
Departments of Neurology; Biostatistics, Bioinformatics & Biomathematics; & Rehabilitation Medicine, Georgetown University, Suite 207 Building D, 4000 Reservoir Road NW, Washington, DC USA 20057

Andrew C. Lee, https://orcid.org/0000-0003-3375-5955
Department of Mathematical Sciences
United States Military Academy
646 Swift Rd
West Point, NY 10996

Rachelle DeCoste, ORCID: 0009-0001-6151-9055
Department of Mathematical Sciences
United States Military Academy
646 Swift Rd
West Point, NY 10996

**Correspondence to:**
Rochelle E. Tractenberg
Georgetown University Neurology
Room 207, Building D
4000 Reservoir Rd., NW
Washington, DC, 20057 USA
**Email: rochelle -dot- tractenberg -at -gmail -dot- com**

Acknowledgement: There are no actual or potential conflicts of interest. Opinions expressed in this article are the authors' own.

Running Head: Mathematical Problem-Solving Pipeline for course and curriculum design

Cite this preprint: Tractenberg RE, Lee AC, Decoste R. (2024). *A Mathematical Problem-Solving Pipeline (MPSP) to strengthen scaffolding in higher education STEM courses*. MathArXiv, https://doi.org/10.48550/arXiv.2412.00009

Abstract

We outline a new tool that can promote coherence within and across higher education mathematics courses by focusing on problem-solving: the Mathematical Problem-Solving Pipeline or MPSP. The MPSP can be used for teaching mathematics and mathematical reasoning authentically. It can be used across courses to develop problem solving skills that can be generalized beyond the course and topic. It has a specific focus on documentation and communication that lets students leverage skills they grow, and use, in other courses or fields. The MPSP can be used in singleton (service) courses, throughout a curriculum, and/or within a capstone experience.

1. Introduction

Mathematical sciences are a critical component of STEM education, playing a vital role in preparing students for a workforce that increasingly relies on computationally intensive fields such as statistics, artificial intelligence, and mathematics. As these disciplines expand across business and government applications, the demand for effective instruction in higher education mathematics has become more pronounced, reflecting its importance in "addressing major challenges in science, technology, and engineering" (Abaté & Cantone, 2005, p. 2/12).

Over the past two decades, there have been significant calls for reform in mathematics education, focusing on innovative, student-centered methods that emphasize conceptual understanding over rote memorization. For example, Abaté & Cantone (2005) highlighted the benefits of contextualizing learning through real-world examples, a key element of constructivist pedagogy. Building on this, Matthews et al. (2010) and Saxena et al. (2016) advocated for the integration of mathematical modeling as a way to engage students and make mathematics more relevant to their experiences. These approaches collectively aim to foster deeper, more meaningful learning.

Abdulwahed et al. (2012) offered a comprehensive overview of the challenges and trends in higher education mathematics, noting a clear shift from procedural or instrumental learning toward conceptual understanding. This shift, grounded in constructivist principles, aligns with efforts to move away from traditional procedural approaches in favor of more engaging and effective instructional strategies. However, while early critiques questioned the value of procedural learning, recent developments suggest an either-or approach to mathematics instruction may not support the achievement of all learning outcomes. Chinofunga et al. (2024) recently documented the use of problem-solving flowcharts as a tool for supporting critical problem-solving learning objectives in mathematics education. This growing interest points to a resurgence of procedural elements—this time incorporated in a way that aligns with conceptual and authentic learning outcomes, reflecting a broader trend toward procedural learning in recent years. Moreover, it can help individual instructors to feasibly update their courses and teaching if they focus on making changes that support specific learning objectives.

While active student engagement is known to enhance mathematical understanding, many classrooms still rely on traditional lecture-based teaching methods, contributing to student struggles and dropout rates in STEM programs. A study out of the University of Nebraska recently showed that 55% of STEM classrooms were lecture-based (Stains, et al. 2018).  Additionally, studies have suggested that students' struggles with mathematics courses can lead them to drop out of college. In their survey of learners who dropped out of community colleges, Ortagus et al. reported that 25.3% indicated that continuing "required math and science courses that were too difficult" among the factors for dropping out (Ortagus et al. 2021). This 2021 result is not surprising, as it is generally accepted that math anxiety can negatively impact student success in the math classroom and beyond. Akin and Kurbangoglu (2011) showed that self-efficacy and positive attitudes (for school persistence) had a negative association with math anxiety. Providing students with a problem-solving approach that they can deploy independently, regardless of the course topic, might allow students who are disengaged in a course, are underprepared, or who have other sources of anxiety about math performance to succeed in courses that tend to impact their overall college success. Increased problem-solving abilities in mathematics could also translate to better problem solving in an unknown future. Kaitera and Harmoinen (2022) noted that instructors may appreciate the need to strengthen students' problem-solving skills in mathematics courses, but that there is a lack support for instructors to innovate. We have adapted a multi-step pipeline of tasks from statistics and data science (Tractenberg, 2020; 2020A) in order to contribute a tool that any mathematics instructor can utilize across diverse course topics and levels: The **Mathematics Problem-Solving Pipeline** (MPSP). The MPSP offers a structured approach that empowers students to develop independent problem-solving skills, reducing math anxiety and fostering success across a range of mathematics courses.

The MPSP has the following tasks:

1. Identify or accept the task/problem to be solved (i.e., some tasks may not be acceptable, e.g., write a solution for an already-solved problem)

2. Develop hypotheses about solutions (structure, methods)

3. Evaluate hypotheses about solutions (structure, methods) (iterative)

4. Apply positivist condition testing (sensitivity analysis)

5. Articulate and test assumptions for solution/method

6. Document what was done in tasks 1-5 (the process and results)

7. Communicate about the work, with appropriate contextualization (and stakeholder consideration)

The MPSP is adapted from the seven tasks that all practitioners in statistics and data science can and do follow and/or recognize (Tractenberg, 2020/2022a), the Statistics and Data Science (SDS) Pipeline, shown in Table 1 for comparison. The SDS Pipeline has the following tasks:

1. Plan/design

2. Collect/munge/wrangle data

3. Analysis – literal for statistics & data science, "evaluation" for computing

4. Interpretation – always for statistics & data science, never for computing

5. Documentation

6. Report & communicate

7. Work on a team

| Main stages of mathematics problem solving | MPSP Tasks | SDS Pipeline Tasks |
|---|---|---|
| Identification of problem and mathematics strategies that can solve the problem | Identify or Accept the Problem; Develop Hypothesis | Plan/Design; Collect/Munge/Wrangle Data |
| Implementation | Evaluate Hypothesis | Analysis |
| Evaluation and Justification | Analyze Sensitivity; Test Assumptions | Interpretation |
| Communication of solution | Document; Communicate | Documentation; Report and Communicate; Work on a team |

Table 1. General alignment of tasks in the MPSP and SPS Pipeline with stages of problem solving (Chinofunga et al., 2024)

These seven tasks are essential to, and recognizable in, the practice of statistics and data science -even if the practitioner does not do each one. The SDS Pipeline is a useful construct for designing and structuring instruction in statistics and data science (Tractenberg, 2022a, 2022b). This structured approach can be effectively adapted for mathematics courses as the MPSP. By employing a consistent instructional framework across diverse mathematics courses, the MPSP facilitates students' ability to build on prior knowledge while integrating new concepts, skills, and abilities as they progress through the curriculum. Furthermore, adapting the SDS framework for mathematics problem-solving enables instructors teaching both mathematics and statistics or data science courses to implement a unified problem-solving methodology across disciplines. Table 2 shows how the MPSP and SDS tasks align along the main stages of mathematical problem solving described by Chinofunga et al (Chinofunga et al., 2024).

To help recall the MPSP tasks, instructors and students can use the acronym **IDEA-DoC,** which encapsulates the seven steps of the MPSP and emphasizes the fact that documentation and communication are core aspects of mathematical work. Each of the letters stand for the following:

**I** - Identify or accept the problem

**D** - Develop hypothesis

**E** - Evaluate hypothesis

**A** - Analyze sensitivity and Test assumptions

**DoC** - Document and Communicate

2. Contexts for utilizing the MPSP to design and document achievement of learning outcomes

In addition to offering structure for all students to use to solve mathematics problems - and for instructors to identify in solutions in student work - the MPSP can also be useful for organizing, providing structure to, and enabling authenticity in problem solving content in core Mathematics courses and curricula. Examples are outlined below.

**2.1 Math Fundamentals.** A range of studies highlight the importance of structure in Calculus and other core curriculum classes. Sullivan et al. discuss how a particular lesson structure can facilitate problem solving and reasoning for students (Sullivan et al., 2015). Bing & Redish (2009) find that students often apply what they call epistemic framing, where they get "stuck" using a limited set of skills or reasoning, even though they possess the necessary tools to solve the problem more effectively (Bing & Redish, 2009). In particular, proper framing and a systematic approach is crucial for students navigating challenging mathematical scenarios. Research has shown that students who employ structured problem-solving strategies demonstrate improved performance and confidence in mathematics. Bell & Polya's classic work on problem solving heuristics still shapes how we approach problems today. It stresses the importance of understanding the problem, making a plan, following the plan, and checking the solution afterward. (Bell & Polya, 1945). Chinofunga et al. explore how procedural flowcharts can support the development of students' problem-solving skills (Chinofunga et al., 2024). A well-known study by Schoenfeld highlights how important it is for students to reflect on their own thinking and manage their actions when solving math problems (Schoenfeld, 2016). Implementing the MPSP throughout these fundamental courses can ensure that students develop and continue to refine their fundamental skills, promoting a coherence through the curriculum or core courses that might otherwise be dependent on individual instructors' choices of approach, which may potentially vary and lack coordination.

**2.2 Math Foundations.** The MPSP provides a unified and consistent approach to implementing instructional innovation targeting problem solving. For example, it facilitates instructors who wish to leverage many diverse options. For example, Reid and Knipping discuss the importance of understanding argumentation structures in mathematics lessons to improve teaching proof processes (Reid and Knipping, 2010). Epp also discusses the challenges students face in writing proofs in mathematics courses, emphasizes the importance of developing logical reasoning skills, and suggests explicit instruction in logic to enhance students' mathematical reasoning abilities (Epp, 2003). Minggi and Mulbar propose a framework for a Local Instruction Theory (LIT) to improve the learning trajectory of students using four broad steps in proof writing: (1) understanding the statement which will be proved, (2) choosing the type of proof, (3) writing the details of the proof, and (4) verifying the validity of the proof (Minggi & Mulbar, 2019). The MPSP brings this structure to all types of mathematics work. Zazkis et al. examined proof-writing behaviors of successful mathematics majors and identified two broad approaches which they termed the targeted strategy and shotgun strategy. When using a targeted approach, students thoroughly understand the statement, select a plan based on this understanding, construct a graphical argument to support the statement's validity, and then formalize this argument into a proof. Conversely, with a shotgun approach, students swiftly experiment with various proof plans and promptly discard a plan at the first sign of difficulty (Zazkis et al., 2015). Integrating the MPSP can bring structure to students' use of either strategy, and instructors can use an MPSP based worksheet to promote evaluable student reflection on which strategy they used and why. Use of the MPSP in courses that are prerequisite for proof writing can prepare students to engage in Minggbi and Mulbar's four steps and foster a focused targeted approach to proofs.

**2.3 Math Modeling**. The process of mathematical modeling involves several key steps that map easily onto the MPSP, including problem identification, hypothesis testing, sensitivity analysis, and communication (Cole et al., 2020; De Corte et al., 2000; Toews, 2012; Zeytun et al., 2017). However, students often struggle with these steps, particularly in recognizing relevant parameters, representing physical situations in equations, and stating justifiable assumptions (Cole et al., 2020). If the MPSP is introduced earlier in the curriculum, then modeling courses can focus more on new steps or methods. A focus on real-world applications, scientific computing, data analysis, and communication skills can enhance the effectiveness of a modeling course (Toews, 2012). Furthermore, the use of open-ended problems and group work can help students develop a deeper understanding of mathematical concepts and improve their problem-solving skills (Greeno et al., 2000). For students who move from fundamentals to foundations to modeling within a curriculum that features the MPSP in multiple contexts, making links from earlier courses to the execution of the same MPSP tasks but with different information and tools can capitalize on prior learning. The MPSP focus on documentation and communication can encourage students to build or utilize their developing communication skills gained in other coursework, reinforcing the "real world-ness" of the modeling.

**2.4 Singleton courses and Capstone Courses.** Many undergraduate degree programs will require only one Mathematics course (or a minimum number), and some students choose these programs specifically due to anxiety they feel about mathematics content. By contrast, students in Mathematics degree programs will often complete a capstone project. In both of these opposing contexts, the incorporation of a recognizable structure like the MPSP can provide support and coherence. The effort needed to integrate the MPSP into any one or set of foundational or service courses would not be wasted, because it can benefit the entire spectrum of students- from anxious students taking their only required Mathematics course to those in the major. Building on a foundation that features the MPSP, the documentation and communication aspects can be featured more prominently in a Capstone context to reinforce skills specifically related to clear communication with non-Mathematics majors. Similarly, documentation and communication skills that non-majors develop from their other coursework can be leveraged to increase the likelihood of catalytic learning (Tractenberg, 2022c) of mathematical argumentation and work.

The seven tasks in the Mathematical Problem-Solving Pipeline are not articulated to imply that every problem will be addressed or solved utilizing all the steps. The purpose of the MPSP is to create a coherent and reproducible structure that instructors can use to direct students' attention to a concrete set of steps they can take – or revisit – whenever they encounter a problem that might be solved with mathematical techniques. The MPSP might help elucidate the "mechanics of mathematics" for all students, whereas in some cases or classes, these are only available to students with the correct intuition about mathematical problem solving. These seven tasks (formulated in 2023) match, and elaborate, the "stages of mathematics problem solving" outlined and studied by Chinofunga et al. (Table 1). Like the SDS Pipeline, the MPSP is consistent with the scientific method, but is more specific (than the SDS Pipeline is) to problem solving. "Work on a team" is a key part of the SDS Pipeline because most applied work in SDS domains involve multi- or cross-disciplinary teamwork. The MPSP was developed principally to support student reasoning and problem solving, and to help all students in STEM courses to develop "mathematical" habits of mind -or, to promote linkage of other knowledge with new mathematical knowledge. The MPSP tasks will be important whenever these habits of mind can be brought to bear. When STEM undergraduates engage in group work, or capstones that involve cross-disciplinary collaboration, or possibly research and internship opportunities, then the MPSP tasks of documentation (for other mathematics practitioners) and communication (of the work for diverse audiences and stakeholders) can support "work on a team" the way it was intended in the SDS Pipeline.

Incorporating the seven MPSP tasks into a STEM course can also facilitate the instructor's achievement of the seven principles of learning in higher education articulated by Ambrose et al. (Ambrose et al., 2010). A Degrees of Freedom Analysis (Tractenberg, 2023) highlights elements that emerge when the principles of learning are aligned with the tasks of the MPSP (Table 2).

Table 2. Alignment of MPSP (columns) with principles of learning (rows; adapted from Ambrose, et al. 2010).

| **MPSP Task:**<br><br>**Principle of Learning:** | Identify the problem | Develop Hypothesis | Evaluate Hypothesis | Sensitivity Analysis | Test Assumptions | Document | Communicate |
|---|---|---|---|---|---|---|---|
| Prior knowledge can be helpful | This task encourages students to review prior work for shared properties to identify structural features of the problem at hand. | Previously encountered hypotheses can be re-considered, re-purposed. | Features of the formulated hypothesis can be leveraged to guide reasoning about how to test the hypothesis, and why testing is needed. | Encourage students to ensure that assumptions are met, method is appropriate, and approximations are plausible - and to consider the impact on the problem/solution **if they are not. *** | Students should be able to confirm that assumptions are plausible and also that they are being met - and advanced students will begin to perceive their **obligations to reconsider the solution if assumptions are not reasonable or met. *** | The majority of mathematics course material will be new (new mathematics content); documenting the technical details of what was done in the previous tasks can help students recognize their own progress and gains in their mathematical reasoning. | Linking students' capabilities in communicating by writing to their growing mathematics abilities can strengthen both. **Responsibilities to communicate clearly and effectively can be considered. *** |
| Knowledge organization supports learning and application of new knowledge | The MPSP offers structure for how mathematics knowledge can be organized and deployed across diverse (authentic) settings and problems. When more than one course utilizes the MPSP, students will develop sophistication for each of the MPSP tasks and will be able to see the best, easiest, and most difficult tasks of the Pipeline from the perspective of each course (e.g., communication about the work in a proof writing course may seem easier than for a calculus course). Utilizing the MPSP to structure problem solving instruction can facilitate each course instructor's integration of real-world problems and problem-based learning into each course in a level- and topically- appropriate way. Each task in the MPSP can be scaffolded over time throughout a program, with decreasing support and increasing independence of the learner over time. ||||||||
| Promotes motivation to learn/sustain learning | Students with experience working through the MPSP tasks can become better self-assessors, diagnosing where in the MPSP they are stronger and weaker. The structure of the Pipeline can also support the formation of a schema for problem solving that can create a catalytic learning experience for students who do not spontaneously develop the motivation to progress (see Tractenberg 2022c). ||||||||

| | |
|---|---|
| Mastery is supported (opportunities to acquire component skills, practice integration, and learn when to apply them) | In courses that feature the MPSP, instructors can begin to share responsibility for student development of a pre-specified target level of performance of each task on the Pipeline. Students can begin to self-assess and also self-direct to finding opportunities to identify and then remediate weaknesses in their performance of each task along the Pipeline. Using the MPSP as a checklist, students can begin to see "a whole solution", e.g., what does it mean when an assessment includes the prompt, "show all work?" This can streamline grading as well - including permitting a coherent peer-evaluation paradigm. |
| Goal-directed practice with formative feedback provided | Using the structure of the MPSP provides sub-goals in the problem solving that is specific to each specific mathematics course topic. Linkages between courses can be made more explicit (e.g., "in single variable calculus you had to <do x> in order to identify the problem, but in multivariable calculus, you need to <do x as well as y> to achieve the same task"). |
| Course climate supports learning | When introducing the MPSP in the course, instructors have a new opportunity to discuss the differences between novice and expert type mathematics practice and problem solving. The structure can also offer an opportunity to discuss scaffolding - by telling students you structure your solutions (utilizing the MPSP) early on, but as the course or program continues, students begin to internalize the steps and self-assess as they complete any problem solution. Considerations of ethical obligations to ensure that a solution is appropriate, that assumptions and approximations (as needed) are supported, and that communication about mathematics work treats all stakeholders- including the discipline of mathematics itself, and other practitioners - fairly. In these ways, the MPSP supports the learning the course was originally designed to promote, while also creating opportunities for students to learn other critical and authentic skills like self-assessment, ethical mathematical practice, and communication about mathematical work. |
| Students will learn to monitor and adjust their approaches to learning. | By breaking down the steps that are general for solving mathematics problems, student attention can be directed to both their capabilities and their own assessment of their capabilities. Instructors can use the MPSP to structure student responses to problems in two ways. First, students fill in a table with their answer for each of the MPSP tasks. Then, students list or evaluate (preferably both) the extent to which their completion of each task meets a priori definitions of "correct" or "competent" performance. |
| **NOTES.** * In these three tasks - and in any course that utilizes the MPSP - instructors can add opportunities for students to consider their work in a social context. Ethical considerations like those indicated (bold) in the table can be integrated into any mathematical course featuring the *d task by, for example, focusing on the impact of the solution, those to whom the solution is communicated, and other stakeholders in the case(s) where the application is not appropriate; the assumptions do not hold; and/or the approximations are not plausible (see, e.g., Tractenberg 2024). | |

Table 2 reflects how courses and instruction that features structure like the MPSP can bring courses into more strong alignment with these key principles of learning. As Bramanti et al. suggests, one way to truly understand a proof or an argument is when we can explain it to somebody else in a conversational style (Bramanti & Travaglini, 2018). However, learning how to communicate - and its importance to the field as well as to other stakeholders - should begin with earlier courses than proof writing. If the MPSP is utilized to structure earlier courses, both documentation and communication, together with the ethical mathematical practice features relating to these important and sometimes underappreciated dimensions of mathematical work, can be integrated and developed across the curriculum. Students who take courses like calculus as part of their core, but who do not major in mathematics, will also benefit from exposure to the MPSP structure because it can prepare them to more fully integrate the mathematics knowledge, skills, and abilities from MPSP-featured courses into their later discipline-specific courses (e.g., Biology).

Table 3 identifies potential assessment opportunities - which in turn suggest instructional ones - featuring each of the MPSP tasks across a variety of mathematics courses typical of both the undergraduate major and service courses that other programs require their students to complete. Table 3 presents an outline of observable behaviors that mathematics students could be expected to perform: a) within each of the identified courses; and b) to demonstrate their capabilities to complete each of the MPSP tasks. As was suggested in Table 1, there are important opportunities to integrate ethical thinking into mathematics courses by attending to the implications of failures in sensitivity analysis, assumptions testing, and communicating work that might remain obscured without the MPSP structure.

Table 3. Example activities within courses for each of the MPSP tasks

| **MPSP Task:** <br> **Course:** | Identify the problem | Develop Hypothesis | Evaluate Hypothesis | Sensitivity Analysis | Test Assumptions | Document | Communicate |
|---|---|---|---|---|---|---|---|
| Single Variable Calculus (Derivatives) | Find the limit of a function as x approaches a value (indeterminate form). Review similar problems | We can substitute the limit value directly | Substitution leads to a form of indeterminacy | Consider the behavior of the function around the x value given an indeterminate form | Test the assumption that the function behaves similarly from both sides of the point | Document approach to simplify the function first before calculating limits | Students can explain a solution in conversational style and can review their peers' explanations critically. Some courses may also elect to have students record their presentations to aid in reflection and evaluation |
| Single Variable Calculus (Integration) | Solving an Integration by parts that loops. Review integration by parts. | Integration by parts formulation | Using the acronym LIATE (Log, Inverse Trig, Algebraic, Trig, Exponential) to set u | Test different values for u and dv | Test conditions that both functions are differentiable | Document approach to solve these types of problems | |
| Multivariable Calculus | Find the dimensions of a rectangular box such that its volume is maximized when inscribed within a sphere of radius r. Review similar problems. | Max volume occurs when dimensions of rectangle form a cube | Maximize volume V=xyz subject to the constraint of a sphere with radius r | How changes in one dimension affects the volume | Validate assumptions about the relationships between dimensions and volume (how to prove that max volume is a cube?) | Document the process step by step, applying proper calculus techniques, with each step clearly shown | |
| Foundations of Math | Proof of Irrationality of √2. Review proofs by contradiction. | If √2 were rational, it could be expressed as the ratio of two integers (p/q), where p and q are coprime integers | Assume √2 is rational and solve | Test other proof approaches | Ensure validity of proof by contradiction Validate assumptions about properties of integers | Document proof step by step, perhaps using a typesetting system like LaTeX to formalize | |
| Math Modeling | Predator-Prey dynamics. Review similar problems. | Hypothesize that it can be described using Lotka-Volterra equations | Data collection, model development, simulate over time | Parameter sensitivity with different environmental scenarios, predation rate, reproduction rate, | Validate assumptions, compare simulated population dynamics with real-world observations, assess reliability and accuracy of data | Document entire modeling process from defining the problem and variables to the model outputs, perhaps using LaTeX to formalize | Present, peer review, or publish findings, tailoring it to the audience |

| | | | | carrying capacity, initial conditions | | | |
|---|---|---|---|---|---|---|---|
| Numerical Analysis | Finding the root of a nonlinear equation. Review similar problems. | Hypothesize that the Newton-Raphson method can approximate the root of a given nonlinear equation | Choose an initial guess and use the Newton-Raphson method and verify convergence | Analyze the sensitivity to the step size in the iterative process and the initial guess | Verify that the function is differentiable in the interval of interest | Document entire method to include problem formulation, method, results, and conclusion, perhaps using LaTeX to formalize | |

As an in-depth example, consider an introductory problem on integration by parts given in a Single Variable Calculus class. Appendices B and C also show examples of a multivariable calculus course and math modeling course. Table 4 shows an assignment template that can be adapted for any course or problem, utilizing the MPSP after the first lesson on this technique. In this example and throughout Table 4, an instructor presents students with the integral:

$$\int e^{-x} \sin x \ dx$$

Table 4. Example Single Variable Calculus Problem Application using Template

| 1. Identify the Problem | $$\int e^{-x} \sin x \ dx$$ Earlier in learning, students are given the actual integral; later they can extract it from a narrative problem statement. |
|---|---|
| 2. Develop Hypothesis | The given integral is similar to previous integrals (such as $\int 4x \cos x \ dx$), so students would hypothesize that the integration by parts formula may apply in this case, using the acronym LIATE (Log, Inverse Trig, Algebraic, Trig, Exponential) as the preferred order to select $u$. |
| 3. Evaluate Hypothesis | To evaluate their hypothesis, students will use $u = \sin x$, and follow through with the integration by parts formula $\int u \ dv = uv - \int v \ du$ where $dv = e^{-x} \ dx$. This yields $du = \cos x \ dx$ and $v = -e^{-x}$, resulting in: $$\int e^{-x} \sin x \ dx = -e^{-x} \sin x + \int e^{-x} \cos x \ dx$$ |
| | Recognizing the need for another integration by parts for $\int e^{-x} \cos x \, dx$, students would proceed with $u = \cos x \, dx$ and $dv = e^{-x} \, dx$. This yields $du = -\sin x \, dx$ and $v = -e^{-x}$. Substituting into the integration by parts formula again: $$\int e^{-x} \sin x \, dx = -e^{-x} \sin x - e^{-x} \cos x - \int e^{-x} \sin x \, dx$$ |
| 4. Sensitivity Analysis | At this point, students might recognize a loop in this problem and believe they chose the wrong $u$, so they might utilize sensitivity analysis to try a different $u = e^{-x}$ instead. However, they would find a similar result. |
| 5. Test Assumptions | They could test assumptions to see if their chosen values of $u$ and $dv$ are indeed differentiable and integrable functions. Once they validate and exhaust all approaches, they will hopefully turn back to the problem and notice the recurrence pattern of $\int e^{-x} \sin x \, dx$ on either side of the equation. Letting $I = \int e^{-x} \sin x \, dx$, student would arrive at $2I = -e^{-x} \sin x - e^{-x} \cos x$), resulting in the final answer of: $$I = \frac{1}{2}(-e^{-x} \sin x - e^{-x} \cos x) + C = \frac{-e^{-x}}{2}(\sin x + \cos x) + C$$ |
| 6. Document | For this problem, the relevant proofs or rules could be documented – this is a narrative part of a solution ("first i tried integration by parts but then i used one-two-buckle-my-shoe"). If multiple methods would yield the same solution, that could be mentioned (possibly for extra credit and/or by more sophisticated students) |
| 7. Communicate | Students summarize the "answer" with a narrative explanation and be able to critique their peers' solutions. For equation-only problems, that is sufficient. For more advanced students, a word problem might have a narrative solution – including limitations or potential uncertainties that accrue from, e.g., using an approximate solution instead of a direct one. |

## 3. Discussion

The MPSP can be utilized to integrate structure into existing mathematics and STEM courses in higher education, to bring such courses into stronger alignment with Principles of Learning (Ambrose et al. 2010) while also addressing key considerations of student-centeredness and the promotion of metacognition in mathematics courses, programs, and interdisciplinary co-instruction (e.g., mathematics & biology). Moreover, utilizing the MPSP to structure homework can facilitate grading and even promote peer review. We conceptualize the MPSP as potentially strengthening scaffolding in higher education STEM courses because it provides the literal scaffold for problem solving, and if utilized across courses in a program, then the independence of students solving mathematical problems can be concretely shown to grow over time.

As discussed earlier, the MPSP can be used in a variety of contexts:
**Singleton course** (the only math course students in another discipline will have to take): Adding the MPSP to a singleton course can add structure (and does not need to add much more than that) so that students can more easily see parallels between problem solving with math and problem solving in their (other) discipline. Rather than adding new information to the course, the MPSP can serve as an organizational feature that can help students to see the logic of mathematical thinking and link this structured approach to other problem solving or inquiry-based projects that they have done or will do in the rest of their courses. Math anxiety - if it cannot be completely overcome, can at least be more concretely specified as students can point to the specific tasks they have the most trouble with. Utilizing the MPSP with students who might be susceptible to math anxiety can actually promote "deliberate practice" (Ericksson, 2007), that is, when the instructor pinpoints exactly what the weakness is and focuses practice and instruction on remediating that weakness.

**Sequence** (a core set of math courses that are foundational knowledge for a major): Using the MPSP can facilitate instructors in course sequences to use the same problem-solving structure as each introduces new material. This enables students to map new information onto an existing framework, so they can focus on making new linkages between old and new knowledge. It also can facilitate assessment development - using templates where the cognitive complexity and sophistication of student responses can be clearly pinpointed as the students' knowledge base and experience grows.

**Degree** (minor/major): As students move through a curriculum towards a degree, problems and projects can be easily mapped onto an MPSP rubric where increasing creativity (iterating between tasks 2-3), deeper comprehension (tasks 2-5), and greater fluency about their work (task 6) can be made visible to the learner, the instructor, and the assessor/evaluator. Capstone courses can feature the MPSP and ask students to reflect on their learning over the program on each task.

Further, the Curriculum Guide to Majors in the Mathematical Sciences (CUPM, Mathematical Association of America, 2015) articulates open questions (OQ) relating to "High-impact" learning experiences (OQ#4):

- **Capstone courses.** What makes an effective capstone course? CUPM continues to call for examples of senior-level capstone courses. Examples will ideally include course objectives, syllabi, requirements, and evaluation of course effectiveness.
- **Research and research-like experiences**. Many different strategies exist to offer undergraduates research or research-like experiences. Such offerings are often expensive; we need new models and alternatives that allow departments to offer such experiences at larger scale. What are alternative ways for students to participate? Can team projects with one advisor and several students provide the desired experience? Might industrial or business sponsors rather than full-time faculty direct team projects that explore applied problems?
- **Internships**. Internships are popular among students. How can the mathematical community encourage partners in business, industry, and government to sponsor more and mathematically richer internships? How can information be shared with advisors and students?

Nine learning objectives for a capstone are identified based on the Boyer Commission Report (1998) and the Educational Effectiveness Working Groups at UC Berkeley (2003). These objectives, which are presented here in a general format so as to be applicable for end-of-degree; end-of-term; and end-of-course capstones, are to:

1. Teach research skills

2. Assess possession of research skills

3. Assess learning of research skills

4. Provide experience with inquiry

5. Assess/estimate independence in research skills

6. Engage inquiry-based learning

7. Teach inquiry-based writing; and that

8A. The Capstone functions formatively; some may also *or instead* be interested in

8B. The Capstone functions summatively.

The MPSP could be used to ensure that mathematics capstone courses function as the Berkeley group envisioned, since the tasks along the Pipeline match scientific inquiry steps (e.g., Wild & Pfannkuch, 1999) including documentation and communication, and are explicitly inquiry based (capstone objectives 4, 6 and 7). As students advance through a degree program, the dependence of MPSP tasks 1-3 on original and creative research skills would then also bring objectives 1-3 into the curriculum so that these could be targeted in the capstone course. If students are asked to fill out a blank table with the MPSP tasks as rows, they could be asked (as mentioned previously) to self-assess their independence (capstone objective 5) and possession (capstone objective 2) of mathematics-specific research skills or of how they are able to harness their mathematical knowledge, skills, and abilities to accomplish research tasks in other domains (e.g, biology). These research-based and structural attributes could be useful in addressing the CUPM (Mathematical Association of America, 2015) open question #4: "High-impact" learning experiences sub-questions.

Luttenberger et al. (2018) discussed the role of attentional systems in the evolution of math anxiety and how they can lead to disruptions in otherwise normally functioning reasoning: "two attentional systems: a top–down, goal-driven system that is influenced by current goals and expectations, and a stimulus-driven system that is influenced by the salient stimuli of the environment." (p. 314). The MPSP could be useful in addressing math anxiety by providing some goal-directed structure for students who find it difficult to initiate problem solving in mathematics courses.

Finally, Abdulwahed et al. (2012) discuss the prominence of a constructivist approach to mathematics instruction. Key features of this pedagogical orientation are:

1. Learning is a student-centered process;
2. Students' autonomy should be fostered;
3. Learning should be contextualized and associated with authentic real-world environments and examples;
4. Social interaction and discourse form an important part of learning;
5. The taught elements should be made relevant to the learner;
6. The taught elements should be linked with the learners' previous knowledge;
7. It is important to facilitate continuous formative assessment mechanisms, self-esteem and motivation;

8. Teachers should act as orchestra synchronizers rather than speech givers; and

9. Teachers should consider multiple representations of their teachings. (Abdulwahed et al. 2012, p. 50).

This list maps easily onto cognitive psychological structural features identified as critical to catalytic learning (Tractenberg, 2022c) and effective teaching by Ambrose et al. (2010) shown in Table 1. Utilizing the MPSP for structure, particularly for multiple mathematics courses taken in sequence, can facilitate instructors' adoption of a constructivist orientation by providing structure within which students can demonstrate their autonomy and their development.

Appendix A. Example Multivariable Calculus Problem Application using Template

| | |
|---|---|
| 1. Identify the Problem | Given a word problem, students must find the dimensions of a rectangular box such that its volume is maximized when inscribed within a sphere of radius $r$. Students can draw a diagram of a sphere with a rectangular box inside it to visualize the problem and help understand the geometric relationships and constraints. |
| 2. Develop Hypothesis | After reviewing similar problems, students should hypothesize that maximum volume of the rectangular box occurs when the dimensions of rectangle form a cube. This simplifies the problem so that the volume of the rectangular box is $V = x^3$. |
| 3. Evaluate Hypothesis | Students would follow through with their hypothesis and note that the diagonal of the box should be equal to the diameter of the sphere, $2r$. Therefore, $\sqrt{3x^2} = 2r$, which results in $x = \frac{2r}{\sqrt{3}} = \frac{2r\sqrt{3}}{3}$. So the volume of the rectangle is $V = x^3 = (\frac{2r\sqrt{3}}{3})^3 = \frac{8r^3\sqrt{3}}{9}$ |
| 4. Sensitivity Analysis | To verify that this is the maximum volume, students can relax the assumption made to maximize the volume $V = xyz$, subject to the constraint $x^2 + y^2 + z^2 = 4r^2$. This problem be solved to optimality using Lagrange multipliers. Doing so would yield the same solution as above. |
| 5. Test Assumptions | Students can test for the assumption in 2-dimensions first to build intuition, or by conducting a numerical or analytical approach to confirm that $V = x^3$ provides the maximum volume. |
| 6. Document | For this problem, students can submit both written along with hand-drawn or technology-aided figures. Documentation can include the problem statement and initial diagram, the hypothesis and reasoning behind it, the mathematical derivation of the solution, and the verification steps (including the Lagrange multiplier method if used). |
| 7. Communicate | Students summarize the answer with a one-sentence narrative explanation or a narrative solution to explain their thought process. This can include a discussion of how this problem relates to real-world applications (e.g., packaging design, efficient use of space in spherical containers) and a reflection on the problem solving process and what was learned. |

Appendix B. Example Math Modeling Problem Application using Template

| 1. Identify the Problem | Given a word problem describing the population dynamics of wolves (predators) and rabbits (prey) in a specific ecosystem, students must develop a mathematical model to predict population changes over time. Students may sketch out the basic relationship between predators and prey to visualize the problem. |
|---|---|
| 2. Develop Hypothesis | Students hypothesize that the Lotka-Volterra equations would be appropriate for modeling this predator-prey system. The basic form of these equations is: $$dx/dt = \alpha x - \beta xy$$ $$dy/dt = \delta xy - \gamma y$$ Where $x$ is the prey population, $y$ is the predator population, and $\alpha, \beta, \delta$, and $\gamma$ are parameters representing birth, predation, predator efficiency, and predator death rates, respectively. |
| 3. Evaluate Hypothesis | Students would set up the Lotka-Volterra equations using the given information from the word problem. They would estimate initial values for the parameters based on the problem description or additional research on wolf-rabbit ecosystems. Using computational tools like Mathematical or Python, they would solve the differential equations numerically and plot the population dynamics over time. |
| 4. Sensitivity Analysis | Students would perform sensitivity analysis by varying the parameters ($\alpha, \beta, \delta, \gamma$) within reasonable ranges to observe how changes affect the model's behavior. This helps identify which parameters have the most significant impact on population dynamics and assess the model's robustness. |
| 5. Test Assumptions | Students can test these assumptions by comparing model predictions with real-world data or by incorporating additional factors to see if they significantly improve the model's accuracy. |
| 6. Document | Students would document their entire process, including the problem statement, hypothesis, mathematical formulation, parameter estimation, numerical solutions, and sensitivity analysis results. They should also include a literature review to ensure their approach is novel or, if not, properly cite existing work on Lotka-Volterra models applied to wolf-rabbit ecosystems. |
| 7. Communicate | Students would prepare a comprehensive report or presentation that includes a clear explanation of the problem and its ecological significance, the mathematical model and its justification, graphical representations of population dynamics, interpretation of results in the context of the ecosystem, discussion of model limitations and potential improvements, and implications for wildlife management or conservation efforts. They should consider stakeholders (e.g., ecologists, wildlife managers, policymakers) and ethical considerations when presenting their findings. |